# Desargues and the "trait à preuves"

Nicolas Bouleau

**Summary**. The 'trait à preuves' is a 17th century expression that refers to drawings rigorously justified by mathematics.

This work is a reflection on the principle of purity of methods from a historical, mathematical and epistemological point of view. Our starting point is the problematic that arose in the 17th century around a dispute between stonemasons and geometers who were theorists of stereotomy. Girard Desargues was involved in this controversy, and we can assume that Desargues' theorem was a key element in his argument. It clearly raises the question of the purity of methods. We will illustrate this with examples. And this leads us to David Hilbert's fundamental work on the subject.

**Résumé**. Le "trait à preuves" est une expression du 17ème siècle qui désigne les tracés rigoureusement justifiés par les mathématiques.

Nous menons dans cet article est une réflexion sur le *principe de pureté des méthodes* d'un point de vue historique, mathématique et épistémologique. Nous partons de la problématique qui s'est nouée au 17ème siècle par une querelle entre les praticiens de la taille de pierre et les géomètres théoriciens de la stéréotomie. Girard Desargues fut impliqué dans cette controverse et on peut penser que le théorème de Desargues constitue pour son auteur un élément clé de son argumentation. Or il pose clairement la question de la pureté des méthodes. Ce que nous illustrons par des exemples. Et cela nous conduit aux travaux fondamentaux de David Hilbert sur le sujet.

## I. Historical background

"*Beauty will result from the form and the correspondence of the whole to the parts, of the parts to each other, and of these to the whole, so that the edifice appears as a whole and well-finished body in which each member is suited to the others and where all the members are necessary to what one wanted to do.*"[1] This recommendation by Andrea Palladio fits in well with the natural philosophy of the relationship between mathematics and architecture during the long period from Achaemenid (Persepolis) and then Greek and Latin Antiquity to the Renaissance, when geometry and arithmetic were available to promote harmony and organise the architect's project. Arithmetic conveys the beauty of music through the simplicity of musical ratios (Alberti).[2] The golden ratio, *divina proportione,* is a trace of what remains divine in human creation (cf. Piero della Francesca,[3] Luca Pacioli[4], Sebastiano Serlio[5]).

Apart from the layout of decorative motifs (metopes and triglyphs, spirals of Ionic capitals, etc.), the basic form is essentially the circle and correlatively the cylinder and the sphere, in Byzantine and Romanesque art, but also for Gothic ogives. The circularity of the arches made it possible, for several labourers to work separately on the voussoirs intended for the same arch.

---

[1] Andrea Palladio *The Four Books of Architecture* (1570) [21]
[2] Leon Baptista Alberti *De re aedificatoria* (1450), Florence 1485 [1]
[3] Piero della Francesca *Trattato d'abaco* 1470 [9]
[4] Luca Pacioli *Summa de arithmetica, geometria, de proportioni et de proportionalita* 1494 [20]
[5] S. Frommel, *Sebastiano Serlio, architecte de la Renaissance,* Gallimard 2002 [10]

But through a series of bold moves that would gradually lead to the Baroque and even the Rococo and Mannerism, the relationship between mathematics and architecture would be reversed in the 16th century in Italy and in the 17th century in the rest of Europe.[6] The taste for complex, unusual shapes, elliptical vaults, oval domes (Francesco Boromini), spiral staircases (inspired by the screw of Saint Gilles), sloping staircases and their connections, ceilings that are feats of construction (Arles town hall by Hardouin-Mansart), called on mathematics as a cry for help that spurred geometers on to develop the theory.[7] This is the birth of a new relationship between mathematics and stonemasons: stereotomy.

This minor art happens to have an important social role, linked to that of perspective, which can make it rank among the *symbolic forms* in the sense of Ernst Cassirer.[8]

## II. The problem of the line and its justification

The arrival of rigorous geometry based on difficult mathematics was not easily accepted by practitioners. A serious controversy arose between the stonemasons on the one hand and the theoretical geometers on the other, which escalated to the point of derogatory invective and the setting up of a 'committee of experts' to decide between the stone cutters' clan, represented by Jacques Curabelle, and the geometers in the person of Girard Desargues.[9]

Desargues' treatise is representative of a perfectly deductive approach in which truly mathematical contributions appear: the study of transformations linked to perspective and quantities *that are invariant* by this type of deformation thanks to the notion of involution and cross ratio. Also the discovery of new properties of conics demonstrated by reducing them to circles (in the case of Blaise Pascal, an emulator of Desargues).[10]

As many historians have noted, and this cannot escape the reader, Desargues' treatise seems truly written to keep away users who are incapable of mathematics.[11] Desargues introduces botanical terminology to evoke his geometric constructions, which is completely useless and could only annoy the stonemasons. It is possible that Desargues, by evoking trees and branches, etc., intended to produce a didactic work, as the encyclopaedists would do in the following century. In this case it was an absolute failure. On the other hand, it was in the field of mathematics that he left a fruitful legacy, projective geometry which was enriched by algebra thanks to complex numbers in Poncelet and Chasles in the 19th century.

This controversy, between a mathematically rigorous approach and a trade, a profession of stonemason knowing how to use the available instruments (as often shown by the complexity

---

[6] A. Chastel *Le grand atelier,* Gallimard 1965 [5]
[7] The example of skewed slopes is followed through the works of various 17th-century architects in the remarkable book of Joël Sakarovitch *Epures d'architectures, de la coupe des pierres à la géométrie descriptive,* Birkhaüser 1998 [26]
[8] Cf. E. Panovsky *Perspective as a symbolic form* (1961) [22] as well as R. Wittkover *The migration of symbols* [31], cf. also J. Lassègue '*Forme symbolique*', [18].
[9] See chapter II section 'L'enjeu de la géométrie' of the book by J. Sakarovitch *op. cit*. p.179 *et seq*.
[10] *La pratique du trait à preuves* par M. Desargues Lyonnais pour la coupe des pierres en architecture, by A. Bosse Paris 1643 [2]
[11] And in this respect differs from earlier works in which perspective is intuitively illustrated, sometimes erroneously, such as H. Vredeman de Vries' *Perspective (*1604) [29]

of the signs of task masons) and knowing that a flat surface can always fit another flat surface, even if it means polishing it afterwards to the desired shape, this dispute necessarily focused on the fact that the designs prior to realisation are always in two dimensions (even if we know that the juries sometimes judged on wood models as for Brunelleschi in the competition for the dome of Florence), whereas it is a question of describing a three-dimensional body and realising it. Basically it is the problem of the map and the territory.

Desargues clearly realised the crucial nature of this transition from two to three dimensions and vice versa. To this end, he included an emblematic statement and its demonstration at the end of his treatise: this is what is commonly known as *Desargues' theorem* (Desargues' other important theorem being his involution theorem[12]).

### III. Desargues' theorem

The statement involves only lines and their intersections. We take two triangles ABC and A'B'C' in the plane such that the lines AA', BB', and CC' are concurrent. The theorem states that in this case the lines bearing the corresponding sides AB and A'B', BC and B'C', CA and C'A' intersect at aligned points. Desargues gives two demonstrations of this. Firstly, by thinking of the figure as the view of an object in space, ABC being the base of a tetrahedron that is cut by a plane according to the triangle A'B'C'. It is clear then that in space the lines AB and A'B' are coplanar and therefore necessarily intersect at the intersection of the planes of the triangles ABC and A'B'C'. And the same goes for the lines BC and B'C' and also CA and C'A'.

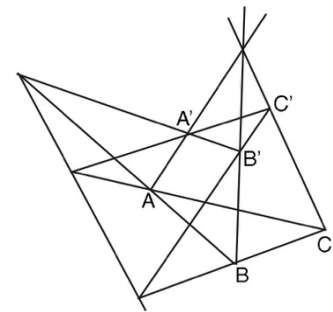
*Desargues' theorem*

Desargues then provides a proof by reasoning in the plane of the initial figure using involution arguments and metric considerations using bi-ratios. The proof by plane geometry reasoning is not immediate. On the other hand, in space, without changing the figure, the property becomes obvious without any calculation.

Curiously, Desargues expresses himself as if the proof in space were not a sufficient argument for the planar result. He gives a detailed proof in two dimensions and insists that this proof also gives a result in three dimensions in whatever way one imagines the planes between the drawn lines.

He takes a view where what is desirable to say, the useful approach to mention, is not at all what we feel today.

Today it is the brevity, the effectiveness and the simplicity of the proof in space that fascinate us and we tend to look for a similar logical phenomenon in other circumstances. For Desargues, one has the impression that he is proud to have found a situation where two-dimensional reasoning 'contains' three-dimensional reality without loss of information. In other words, he highlights a case where the plan view of the structure corresponds perfectly with reality, which adds weight to his argument with regard to the masons. Mr Poudra, the 19th-century publisher of the treatise, expressed it as follows: "Thus, on the plane of the paper, there results a figure which corresponds point to point, straight line to straight line and cause to cause, to that of the various planes, and then one can discuss their properties on one as on the other and by this means dispense with that of the relief, substituting for it those of a single plane. An important observation that reveals the goal Desargues sets out to achieve in this proposition."[13]

---

[12] Cf. J. Fauvel & J. Gray [8] pp368-370.

Desargues' goal is to prove that plane figures can perfectly master certain non-elementary spatial configurations.

## IV. Two examples

The first describes a situation that generalises the case of Desargues' triangles, bearing in mind that today, as Michel Chasles wrote, thanks to complex projective geometry "*therefore anyone who wants to, in the current state of science, can generalise and create in geometry; genius is no longer indispensable to add a stone to the edifice*".[14]

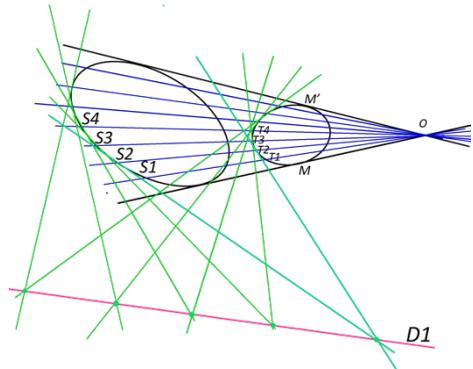

Let there be two ellipses $E_1$ and $E_2$ in the plane. We trace the two common tangents that leave the ellipses on the same side. They intersect at O. This determines two arcs on each ellipse, namely *a* and *a'* for the first and *b* and *b'* for the second.

The proposition states that if we consider several secants originating from O that intersect arc *a* at points $S_n$ and arc *b* at points $T_n$, then the points of intersection of the tangents at $S_n$ and $T_n$ are aligned on a straight line, say $D_1$. And if the same secants cut *a'* at $S'_n$ and *b'* at $T'_n$, the tangents at $S'_n$ and $T'_n$ intersect at points on the same straight line $D_1$. Finally, if we pair *a* with *b'* and *a'* with *b*, the tangents constructed in a similar way intersect on a straight line $D_2$.

*Proof.* The ellipses can be considered to be on a horizontal plane. Consider the cone with base $E_1$, whose vertex is on the vertical line through *O*, and the vertical cylinder with base $E_2$. This cone and this cylinder intersect according to a biquadratic curve that has two double points and is therefore reduced to two plane curves that are ellipses.[15] The pairs of tangents designated in the statement are the projections of coplanar tangents since in a plane tangent to the cone according to a generatrix. In space, these 

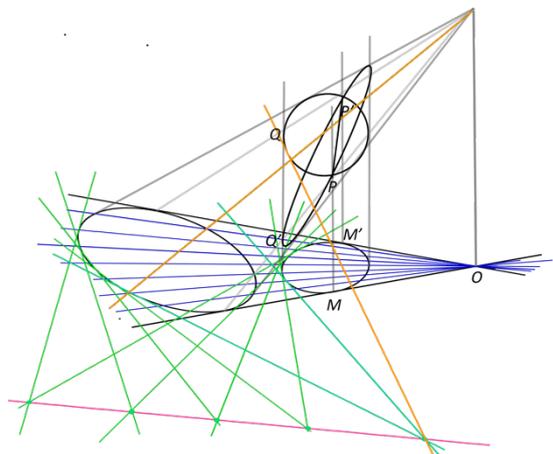

tangents necessarily intersect at the intersection of the horizontal plane and the plane of the ellipse PP'Q.[16]

---

[13] *Œuvres de Desargues, réunies et analysées par M. Poudra* 1864 [25]
[14] Chasles M., *Aperçu historique sur l'origine et le développement des méthodes en géométrie*, Brussels, (1837) [4]
[15] Ernest Duporcq *Premiers principes de géométrie moderne*, [7] p86.
[16] This last argument similarly shows that if a ruled developable surface is cut by two planes, the tangents to the two intersection curves at points on the same generatrix intersect at aligned points when the generatrix varies.

We can replace *a* with *a'* and *b* with *b'* by reasoning again with the plane PP′Q. And we can use *a* and *b'* as well as *a'* and *b* by reasoning with the plane PP′Q′, since the sections by PP′Q and PP′Q′ are projected along the same ellipse E$_2$.

Daniel Perrin told me that the property may be generalised. It extends to two conics where two of the points of intersection are imaginary conjugates by reducing to two circles, and for the case where the conics intersect in four real points, a proof can be established by finding a homography that transforms one conic into the other and fixes the intersection of the common tangents. [17]

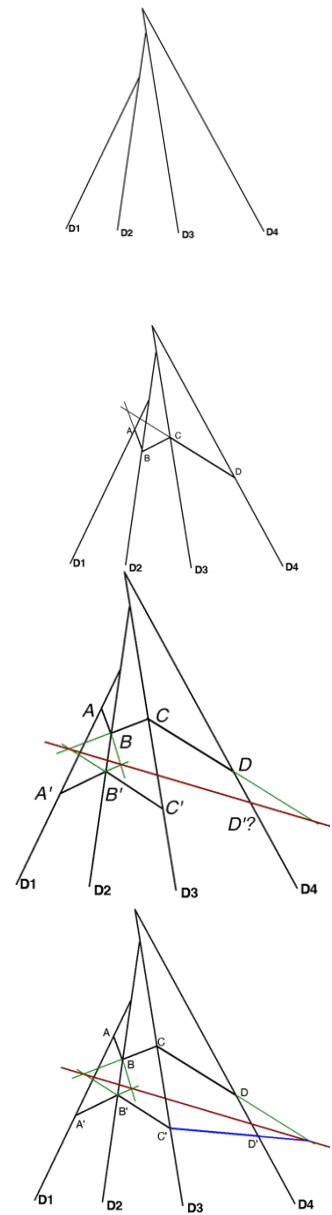

A second example can be presented as follows to simplify the explanation. Consider the surface *S* in space obtained by folding a sheet of paper along lines *d*1, *d*2, *d*3 and *d*4, the faces between these edges being assumed to be flat. Such a surface is a basic example of a developable surface. By central or orthogonal projection (as in descriptive geometry), this surface gives four lines *D*1, *D*2, *D*3 and *D*4. A useful question for the draughtsman is to know under what conditions a quadruplet of points *A, B, C, D*, with *A* on *D*1, *B* on *D*2, *C* on *D*3 and *D* on *D*4, can be considered as the projection of the vertices *a* on *d*1, *b* on *d*2, *c* on *d*3 and *d* on *d*4, of a section of *S* by a plane.

Let us assume that this is the case for the points A, B, C, D. That is to say that the lines *ab* and *cd* intersect in space (or are parallel) and are therefore coplanar.

And this being assumed, let us consider a second quadruplet *A'B'C'D'* with *A'* still on *D*1, *B'* on *D*2, *C'* on *D*3 and *D'* on *D*4. Is it possible to raise back this quadruplet in a planar section of *S* after we have already done so for *ABCD*?

This time the quadruplet *A'B'C'D'* must satisfy a condition: the point *D'* is determined if the points *A', B', C'* are chosen. The situation is unusual and could shock the "common sense" of the draughtsman: if *A, B, C, D* have already been chosen as a projection of a plane section then *D'* is determined by the choice of *A', B', C'*.

The *proof* is based on the argument that if *ABCD* and *A'B'C'D'* are projections of two plane sections of *S* then, because the faces of *S* are plane, the points of intersection of the lines *AB* and *A'B'*, *BC* and *B'C'*, *CD* and *C'D'*, are aligned on the projection of the intersection of the planes of the two plane sections.

These two examples show that reasoning in space provides logical implications that were not available if we remained in plane geometry. Reasoning in space can sometimes lead to results that were not immediately apparent in the restricted framework, a point on which Hilbert shed much light.

---

[17] Cf. https://www.imo.universite-paris-saclay.fr/fr/perso/daniel-perrin/

## V. David Hilbert returns to the task

In mathematics, the 19th century was characterised by several currents of thought that led to a revival of axiomatic thinking, in which geometry played a permanent role. After the appearance of non-Euclidean geometries with Gauss, Bolyai and Lobatchevski, and the use of transformations by Jean-Victor Poncelet and Felix Klein (Erlangen programme), reasoning by translation from one geometry to another appeared, establishing *relative logical* coherence between various choices of primitive notions.

During the last ten years of the century, an abundance of axiomatic systems were defined for geometry and arithmetic. M. Pieri proposed an axiomatisation of projective geometry.[18] In 1891, H. Wiener detached the use of terms (point, line, etc.) from their usual meaning. At the International Congress of Philosophy in Paris in 1900, A. Padoa presented an axiomatisation of arithmetic.

But it was the study conducted by David Hilbert in his *Grundlagen der Geometrie* [19] that made the greatest impression on people at the time. In it, he analyses the axioms of various geometries to determine their dependence or independence. In particular, he shows that the *theorem of Desargues* can be false in certain geometries related to two-dimensional Euclidean geometry (non-Arguesian geometries). He shows that certain congruence properties that can be derived from the third dimension are necessary to obtain the plane theorem of Desargues.

Returning to the problem that arose in the 19th century, we see that Desargues' two-dimensional demonstration uses notions of congruence and equality relations that are not verified by all possible interpretations of the basic terms and their relations. Hilbert thus establishes that for this problem the use of the third dimension is equivalent to adding postulates to plane geometry.

This means, if the axioms are made explicit, that the *principle of purity of methods* is questionable in the case of Desargues' theorem.[20]

The *principle of purity of methods* is generally stated as a methodological rule, a *rule of the art*, a bit like consonance in tonal music. It was considered preferable, more beautiful let's say, more elegant, to use in the proof of a theorem only notions of the same nature as those necessary for the understanding of the statement of the theorem. This methodological recommendation stemmed from the fact that until the end of the 19th century, mathematics was not formalised, so that it was not absolutely clear what mathematics really was; it was limited by somewhat grey areas. Mathematics was considered the language of the great treatises, those of Euclid, Descartes, etc. It was therefore advisable to use hypotheses sparingly so as not to add outside knowledge.

At the very beginning of the 20th century, approaches to axiomatisation multiplied and philosophical biases became more radical.

---

[18] On the emergence of projective geometry cf. J. Fauvel, J. Gray, *The History of Mathematics, a Reader*, Macmillan Press 1987 [8] p366 *et seq*.
[19] David Hilbert <u>Grundlagen der Geometrie (</u>1899) <u>[16]</u>.
[20] A similar hetero-poietic detour demonstrates Monge's theorem of the three circles. Cf. P. Mancosu *Infini, logique, géométrie,* [19] p.352 *et seq.*

Hilbert emphasised the objective of demonstrating the coherence of the axioms used (the second problem stated at the Paris symposium in 1900) and specified this project in 1904, then again at the beginning of the 1920s with his students W. Ackerman, P. Bernays and J. Herbrand, especially since, beyond the quarrels with Brouwer and the supporters of intuitionism, everyone was convinced that some of the intuitionists' concerns were well-founded, and that it was particularly important to rely firmly on finite notions and very cautiously on those that use infinity.

In different forms, the *programme of Hilbert* to demonstrate the coherence of mathematics and already of arithmetic consists of focussing attention on the demonstrations themselves as an object of study (metamathsmatics) in order to construct what can be called a progressive finitism which establishes that the passage from demonstrated formulae to new ones is done by means of things already solidly established and finitist enumerations. As for abstract notions, they should be accepted for convenience, provided that it can be shown that they can be dispensed with by returning to indisputable paths leading to the same results. [21]

In other words, Hilbert's programme is based on a *principle of purity of methods* allowing one to dispense with anything that is not *elementary and finitely* determinate. The Herbrand theorem is an example of a result that would be acceptable in a barely relaxed form of Hilbert's programme.[22]

Why did Hilbert take this route when he knew that impure methods were often more powerful, as he demonstrated in the case of geometry? Because he had no choice. If he had used notions of transfinite analysis or even elementary analysis to establish the consistency of arithmetic, he would have contributed nothing of real value. There is an obvious *bootstrap* problem. Nevertheless, it was not absurd to consider that by studying the assemblages of signs, their construction, and their sequences by logic, that is to say discrete things, metamathematics could contribute more than remaining at the level of mathematics itself, where the meaning of the concepts deceives us about their *combinatorial* complexity based on axioms.

The rest of the History of Ideas, with Gödel, Church and Turing, showed that *combinatory logic* was at the heart of what makes mathematics creative. The failure of Hilbert's programme is a break that opens up a new philosophy of mathematics where combinatory logic plays a major role in 'innovativity'. And this is also instructive in chemistry and molecular biology.[23]

Without going into the subject here, let us make a few remarks. An important observation is that Peano's axioms for arithmetic are extremely convincing. It is not clear what could be improved on them to talk about whole numbers, neither adding nor subtracting. And yet it is an incomplete system. So what can happen?

The set of theorems, i.e. demonstrable statements, should not be considered as given. We do not have explicit access to them. The predicate that defines it after enumeration of the statements is of the form ∃*p* *B*(*p*, *n*) where *p* is the number associated with the proof of the statement with number *n*. It is a recursively enumerable non-recursive predicate. This means that it is beyond the mechanisable. Given statement of number *n*, there is no effective bound as a function of *n* on the length of proofs that must be tried to demonstrate it.

---

[21] For a discussion of the links between Hilbert's programme and the principle of method purity, see Jean-Yves Girard 'Le théorème d'incomplétude de Gödel' *in N.* Bouleau, J.-Y. Girard and A. Louveau *Cinq conférences sur l'indécidabilité* [11]

[22] Herbrand describes a rule of formation associating with each formula A of first-order predicate calculus with quantifiers an infinite countable family of formulas without quantifiers such that A is deductible if and only if a formula of the family is.

[23] Cf. N. Bouleau *Ce que Nature sait,* 540 p., Presses Universitaires de France 2021; also "Scientisme et environnement : l'opportunisme stratifié" *5ème Journée Condorcet*, Montpellier University. 5 nov. 2024, online.

This discovery can be expressed by saying that there is a victory of combinatorics over nomology, that is to say over what is governed by laws (from *νόμος*, the law). There are assemblages that have no other explanation than to be themselves and nevertheless provide shortcuts that change very long proofs into much shorter ones.

Finally, the dispute encountered by Desargues, and Hilbert's courage in designating the path he considered most important despite its difficulty, have made it possible to renew epistemology by suggesting that we place 'normal' science in the sense of Thomas Kuhn as a principle of purity of methods that *revolutions in* science sometimes undermine.